\documentclass[4pt, oneside]{article}
\usepackage{amsmath,amssymb,amsthm,units,stmaryrd}

\usepackage{yfonts}
\usepackage{amsmath,amssymb,amsthm,units,stmaryrd}
\usepackage{graphicx}
\usepackage[all]{xy}
\newtheorem{conj}{Conjecture}[section]
\newtheorem{definition}{Definition}[section]
\newtheorem{lemm}{Lemma}[section]
\newtheorem{prop}{Proposition}[section]
\newtheorem{cor}{Corollary}[section]
\newtheorem{theorem}{Theorem}[section]

\def\lb{\left\llbracket}
\def\rb{\right\rrbracket}

\def\rra{\mapsto}

\def\fmodels{\xymatrix{
\ar@{|=}[r]^{<\omega}&
}
}
\def\nmodels{\xymatrix{
\ar@{|=}[r]^{N}&
}
}
\def\<{\left <}
\def\next{{f}}
\def\hf{{[f]}}
\def\ev{{\left< f\right>}}
\def\nc{{\Box}}
\def\ps{{\Diamond}}

\def\peq{\preccurlyeq}
\def\seq{\succcurlyeq}
\def\>{\right >}
\def\bra {\left [}
\def\ket{\right ]}
\def\cbra{\left \{}
\def\cket{\right \}}

\DeclareSymbolFont{AMSb}{U}{msb}{m}{n}
\DeclareMathSymbol{\N}{\mathbin}{AMSb}{"4E}
\DeclareMathSymbol{\Z}{\mathbin}{AMSb}{"5A}
\DeclareMathSymbol{\R}{\mathbin}{AMSb}{"52}
\DeclareMathSymbol{\Q}{\mathbin}{AMSb}{"51}
\DeclareMathSymbol{\I}{\mathbin}{AMSb}{"49}
\DeclareMathSymbol{\C}{\mathbin}{AMSb}{"43}
\begin{document}
\date{\today}
\title{A sound and complete axiomatization for Dynamic Topological Logic
      }
\author{David Fern\'{a}ndez-Duque
\\ {\small Group for Logic, Language and Information,}
\\ {\small Universidad de Sevilla}$  $
        }
\maketitle
\begin{abstract}
Dynamic Topological Logic ($\mathcal{DTL}$) is a multimodal system for reasoning about dynamical systems. It is defined semantically and, as such, most of the work done in the field has been model-theoretic. In particular, the problem of finding a complete axiomatization for the full language of $\mathcal{DTL}$ over the class of all dynamical systems has proven to be quite elusive.

Here we propose to enrich the language to include a polyadic topological modality, originally introduced by Dawar and Otto in a different context. We then provide a sound axiomatization for $\mathcal{DTL}$ over this extended language, and prove that it is complete. The polyadic modality is used in an essential way in our proof.
\end{abstract}

\section{Introduction}

Dynamic Topological Logic ($\mathcal{DTL}$) is a combination of topological modal logic \cite{tarski} and temporal logic \cite{temporal} used for reasoning about {\em dynamic topological systems}, which are pairs $\<X,f\>$ consisting of a topological space $X$ and a continuous function $f:X\to X$. $\mathcal{DTL}$ was introduced in \cite{arte} as $\mathsf{S4C}$; later \cite{kmints} added an infinitary temporal modality, here written $\hf$ (`henceforth'), into the language. This development allows us to reason about arbitrary iterations of $f$ and capture long-term phenomena such as recurrence.

Although a substantial body of work has been done on the logic, due to its model-theoretic definition, it has proven difficult to work with it in an entirely syntactic manner. An axiomatization has been suggested \cite{kmints}, but it has not been proven complete. However, the logic is recursively enumerable \cite{me2}, which gives hope of finding a proof system for it.\footnote{It is, nevertheless, undecidable \cite{konev}.}

In this paper we present an extension to $\mathcal{DTL}$, which we shall call $\mathcal{DTL}^\ast$, generalizing the use of the topological modal operator to its `tangled' variant, introduced in \cite{do} and also studied in \cite{me:simulability,me:tangle}. Thus we obtain expressions of the form $\ps\Gamma$, where $\Gamma$ is a finite set of formulas; the ordinary monadic modality becomes a special case when $\Gamma$ is a singleton, and we write $\ps\gamma$ for $\ps\cbra\gamma\cket$. The interpretation of these formulas uses the {\em tangled closure} operator, discussed in Section \ref{sectan}.

The axiomatization is mostly the amalgamation of proof systems for the isolated modalities, but we need a polyadic version of the usual continuity axiom from \cite{kmints} which takes the form
\[\ps\next\Gamma\to\next\ps\Gamma.\]

This axiom is sound \cite{dynamictangle}; our main goal is to show that the proof system is complete. For this we expand on techniques from \cite{me2}, where simulations are an essential tool; it is because simulability is expressible with the tangled operator (but not in the ordinary modal language) that we need to use this enriched language. Our completeness proof relies heavily on \cite{me:simulability,dynamictangle}, where we began analyzing the topological behavior of the tangled modality.

In \cite{me:simulability}, we show that given a finite, pointed $\mathsf{S4}$ model $\mathfrak w$, the property of being {\em simulated} by $\mathfrak w$ is not always definable over the class of topological models in the basic modal language; however, using the polyadic modality, there is always a formula\footnote{The language used in \cite{me:simulability} is different, but expressively equivalent, to the purely topological fragment of the language we will use here.} $\mathrm{Sim}(\mathfrak w)$ which defines being simulated by $\mathfrak w$. This is essential to our current completeness proof, since simulations play a key role in \cite{me2} and capturing them syntactiaclly is an important step in our argument.

Meanwhile, \cite{dynamictangle} gives a sound and complete axiomatization for the polyadic $\mathsf{S4C}^\ast$, that is, the fragment of $\mathcal{DTL}^\ast$ without $\hf$. Note that the monadic $\mathsf{S4C}$ was proven complete in \cite{arte}.

While a full axiomatization of (an extension of) $\mathcal{DTL}$ as we are presenting here is novel, there are many positive and negative results regarding axiomatizations of related systems. We summarize them below:
\begin{description}
\item[The next-interior fragment $\mathsf{S4C}$.] This fragment uses only the monadic mo\-da\-lity $\ps$ and the next-time operator $\next$. A sound and complete axiomatization is given in \cite{arte}, where the logic is also shown to be decidable. A complete axiomatization for the logic over spaces with {\em homeomorphisms} is given in \cite{kmints}.
\item[Monadic $\mathcal{DTL}$ over arbitrary systems.] This logic is undecidable \cite{konev} but recursively enumerable \cite{me2}. However, the recursive enumeration does not suggest a reasonable axiomatization. It is conjectured in \cite{kmints} that a rather intuitive proof system is complete; however, this has never been proven.

It should be remarked that the only $\mathcal{DTL}$'s with the `henceforth' modality which have been given a proper axiomatization are over trivial spaces (where the only open sets are the empty set and the entire space) and almost disjoint spaces (where every open set is closed); this can be found in \cite{s5}.
\item[Monadic $\mathcal{DTL}$ over spaces with homeomorphisms.] When restricting semantics to spaces with homeomorphisms but allowing the `henceforth' modality in the language, $\mathcal{DTL}$ becomes non-axiomatizable \cite{wolter}.
\item[The polyadic $\mathsf{S4C}^\ast$.] Here we use only $\ps$ and $\next$, but $\ps$ is allowed to act on finite sets of formulas, and is interpreted as a `tangled closure' operator. New axioms are needed to define the behavior of $\ps\Gamma$ and to describe the interaction of the dynamics with the polyadic modality; they appear in our axiomatization in Section \ref{theax}. This logic is then proven sound and complete in \cite{dynamictangle}.

\end{description}

The layout of this paper is as follows. Section \ref{topre} reviews topologies and their relation to preorders; this relation is useful in linking topological and Kripke semantics of $\mathsf{S4}$. Section \ref{sectan} reviews the {\em tangled closure} operator, which is an important addition to the expressiveness of $\mathcal{DTL}$. Section \ref{basic} gives the formal language and its semantics, and then Section \ref{theax} describes our proposed axiomatization. Subsequent sections mainly review notions from \cite{me2}, although with some modifications to accommodate the new tangled modality: Section \ref{secndq} gives an overview of quasimodels, Section \ref{secgen} shows how one obtains limit models from quasimodels, Section \ref{secsim} discusses simulations and Section \ref{secfour} introduces the universal state space. 

In Section \ref{secsubchar} we discuss the properties of the formulas $\mathrm{Sim}(\mathfrak w)$; this section depends on results from \cite{me:simulability,dynamictangle}. With this we define canonical quasimodels in Section \ref{seccan}. Section \ref{re} reviews efficiency, originally used within the context of $\mathcal{DTL}$ in \cite{konev} and an important tool in showing that it is recursively enumerable in \cite{me2}. In Section \ref{tempinc} we use these ideas to show that canonical structures are in fact quasimodels: this is used in Section \ref{secomp}, where our main completeness result is stated and proved. Finally, Section \ref{conc} gives an outlook for future work and discusses a possible application to the $\mathcal{DTL}$ of minimal systems.

\section{Topologies and preorders}\label{topre}

The purpose of {\em Dynamic Topological Logic} is to reason about dynamical systems defined over topological spaces. Such spaces provide an interpretation of the modal logic $\mathsf{S4}$, generalizing its well-known Kripke semantics.

Let us recall the definition of a {\em topological space}:

\begin{definition}[topological space]
A {\em topological space} is a pair $\mathfrak X=\<|\mathfrak X|,\mathcal{T}_\mathfrak X\>,$ where $|\mathfrak X|$ is a set and $\mathcal T_\mathfrak X$ a family of subsets of $|\mathfrak X|$ satisfying
\begin{enumerate}
\item $\varnothing,|\mathfrak X|\in \mathcal T_\mathfrak X$;
\item if $U,V\in \mathcal T_\mathfrak X$ then $U\cap V\in \mathcal T_\mathfrak X$ and
\item if $\mathcal O\subseteq\mathcal T_\mathfrak X$ then $\bigcup\mathcal O\in\mathcal T_\mathfrak X$.
\end{enumerate}

The elements of $\mathcal T_\mathfrak X$ are called {\em open sets}. Complements of open sets are {\em closed sets}.
\end{definition}

Given a set $A\subseteq |\mathfrak X|$, its {\em interior}, denoted $A^\circ$, is defined by
\[A^\circ=\bigcup\cbra U\in\mathcal T_\mathfrak X:U\subseteq A\cket.\]

Dually, we define the closure $\overline A$ as $|\mathfrak X|\setminus(|\mathfrak X|\setminus A)^\circ$; this is the smallest closed set containing $A$.

Topological spaces generalize transitive, reflexive Kripke frames. Recall that these are pairs $\mathfrak W=\<|\mathfrak W|,\peq_\mathfrak W,\lb\cdot\rb_\mathfrak W\>$ where $\peq_\mathfrak W$ is a preorder on the set $|\mathfrak W|$. We will write $\peq$ instead of $\peq_\mathfrak W$ whenever this does not lead to confusion.

To see a preorder as a special case of a topological space, define
\[\mathop\downarrow w=\cbra v:v\peq w\cket.\]
Then consider the topology $\mathcal T_\peq$ on $|\mathfrak W|$ given by setting $U\subseteq|\mathfrak W|$ to be open if and only if, whenever $w\in U$, we have $\mathop\downarrow w\subseteq U$ (so that all sets of the form $\downarrow w$ provide a basis for $\mathcal T_\peq$). A topology of this form is a {\em preorder topology}\footnote{Or, more specifically, a {\em downset topology}.}. It is not hard to check that the Kripke semantics given by $\peq$ coincide with the topological semantics given by $\mathcal T_\peq$.

Throughout this text we will often identify preorders with their corresponding topologies, and many times do so tacitly.

We will also use the notation
\begin{itemize}
\item $w\prec v$ for $w\peq v$ but $v\not\peq w$ and
\item $w\sim v$ for $w\peq v$ and $v\peq w$.
\end{itemize}

The relation $\sim$ is an equivalence relation; the equivalence class of a point $x\in X$ is usually called a {\em cluster}, and we will denote it by $[x]$.

\section{The tangled closure}\label{sectan}

In this paper we will enrich the language of $\mathcal{DTL}$ by a topological operator called the {\em tangled closure}, which generalizes the ordinary closure to families of sets and not only single sets. It was introduced in \cite{do} for Kripke frames and has also appeared in \cite{me:simulability,me:tangle,dynamictangle}.

\begin{definition}[Tangled closure]
Let $\mathfrak X$ be a topological space and $\mathcal S\subseteq 2^{|\mathfrak X|}$.

Given $E\subseteq|\mathfrak X|$, we say $\mathcal S$ is {\em tangled} in $E$ if, for all $A\in \mathcal S$, $A\cap E$ is dense in $E$.

We define ${\mathcal S}^\ast$ to be the union of all sets $E$ such that $\mathcal S$ is tangled in $E$.
\end{definition}

It is important for us to note that the tangled closure is defined over any topological space; however, we will often be concerned with locally finite preorders in this paper. Here, the tangled closure is relatively simple.

\begin{lemm}
Let $\<S,\peq\>$ be a finite preorder, $x\in S$ and $\mathcal O\subseteq\mathcal P(S)$. Then, $x\in{\mathcal O}^\ast$ if and only if there exist $\<y_O\>_{O\in\mathcal O}$ such that $y_O\peq x$ for all $O\in\mathcal O$ and $y_O\sim y_{O'}$ for all $O,O'\in \mathcal O$.
\end{lemm}

\proof
Suppose $\mathcal O=\cbra O_i\cket_{i<I}$ and $x\in\mathcal O^\ast$.

If $y\in{\mathcal O}^\ast$ and $i<I$, we have that $y\in\overline{O_i\cap{\mathcal O}^\ast}$ (because $\mathcal O$ is tangled in $\mathcal O^\ast$ \cite{dynamictangle}) which means that there is $z\peq y$ with $z\in O_i\cap{\mathcal O}^\ast.$ We can apply this fact countably many times to obtain a sequence
\[x\seq y_0\seq y_1\seq...\]
such that $y_n\in O_i\cap \mathcal O^\ast$ if and only if $n\equiv i\pmod{I}$.

Since $S$ is finite the sequence must eventually stabilize, which means that for some value of $N$ we have that $y_n\sim y_m$ for all $n,m>N$. Thus
\[\cbra y_{N+i}:i<I\cket\]
is a set of points which are all equivalent in $\peq$ and contain a representative of each $O_i$, as desired.
\endproof

\section{Dynamic Topological Logic}\label{basic}

The language of $\mathcal{DTL}^\ast$ (henceforth $\mathsf L^\ast$) is built from propositional variables in a countably infinite set $\mathsf{PV}$ using the Boolean connectives $\wedge$ and $\neg$ (all other connectives are to be defined in terms of these), the unary modal operators $\next$ (`next') and $\hf$ (`henceforth'), along with a polyadic modality $\ps$ which acts on finite sets\footnote{We use $\ps$ as primitive rather than $\nc$ because we find its meaning more intuitive.}, so that if $\Gamma$ is a finite set of formulas then $\ps\Gamma$ is also a formula. Note that this is a modification of the usual language of $\mathcal{DTL}$. We write $\nc$ as a shorthand for $\neg\ps\neg$; similarly, $\ev$ denotes the dual of $\hf$. We also write $\ps\gamma$ instead of $\ps\cbra\gamma\cket$; its meaning is identical to that of the usual $\mathsf{S4}$ modality \cite{dynamictangle}.

We will denote fragments of $\mathsf L^\ast$ by indicating the modalities which are allowed in them; for example, $\mathsf L^\ast_{\ps}$ is the language of polyadic unimodal logic and $\mathsf L^\ast_{\ps f}$ is the fragment without $\hf$ corresponding to $\mathsf{S4C}^\ast$ \cite{dynamictangle}. The star indicates the use of the polyadic $\ps$, so that, for example, $\mathsf L_\ps$ denotes the standard modal language.

Formulas of $\mathsf L^\ast$ are interpreted on dynamical systems over topological spaces, or {\em dynamic topological systems}.

\begin{definition}[dynamic topological system]
A {\em dynamic topological system (dts)} is a triple
\[\mathfrak X=\<|\mathfrak X|,\mathcal{T}_{\mathfrak X},f_{\mathfrak X}\>,\]
where $\<|\mathfrak X|,\mathcal{T}_{\mathfrak X}\>$ is a topological space and
\[f_{\mathfrak X}:|\mathfrak X|\to |\mathfrak X|\]
is continuous.
\end{definition}

\begin{definition}[valuation]
Given a dynamic topological system $\mathfrak X$, a {\em valuation} on $\mathfrak X$ is a function
\[\lb\cdot\rb:\mathsf L^\ast\to 2^{|\mathfrak X|}\]
satisfying
\[
\begin{array}{lcl}
\lb\alpha\wedge\beta\rb_{\mathfrak X}&=&\lb\alpha\rb_{\mathfrak X}\cap \lb\beta\rb_{\mathfrak X}\\\\
\lb\neg\alpha\rb_{\mathfrak X}&=&|\mathfrak X|\setminus \lb\alpha\rb_{\mathfrak X}\\\\
\lb\next\alpha\rb_{\mathfrak X}&=&f^{-1}\lb\alpha\rb_{\mathfrak X}\\\\
\lb\hf\alpha\rb_{\mathfrak X}&=&\displaystyle\bigcap_{n\geq 0}f^{-n}\lb\alpha\rb_{\mathfrak X}\\\\
\lb\ps\cbra\alpha_1,...,\alpha_n\cket\rb_\mathfrak X&=&\cbra\lb\alpha_1\rb_\mathfrak X,...,\lb\alpha_n\rb_\mathfrak X\cket^\ast.\\\\
\end{array}
\]
\end{definition}

We may also write $\lb\Gamma\rb_\mathfrak X$ instead of $\cbra\lb\gamma\rb_\mathfrak X:\gamma\in\Gamma\cket$, so that
\[\lb\ps\Gamma\rb_\mathfrak X=\lb\Gamma\rb^\ast_\mathfrak X.\]

A {\em dynamic topological model} (dtm) is a dynamic topological system $\mathfrak X$ equipped with a valuation $\lb\cdot\rb_\mathfrak X$. We say a formula $\varphi$ is {\em valid} on $\mathfrak X$ if $\lb\varphi\rb_\mathfrak X=|\mathfrak X|$, and write $\mathfrak X\models\varphi$. If a formula $\varphi$ is valid on every dynamic topological model, then we write $\models\varphi$.

We will often write $\<\mathfrak X,x\>\models\varphi$ instead of $x\in\lb\varphi\rb_{\mathfrak X}$.

\section{The axiomatization}\label{theax}

Our proposed axiomatization for $\mathcal{DTL}^\ast$ consists of the following:

\begin{description}
\item[$\mathsf{Taut}$] All propositional tautologies.
\item[]Axioms for $\ps$:
\begin{description}
\item[$\mathsf K$] $\nc(p\to q)\to(\nc p\to\nc q)$
\item[$\mathsf{T}$] $\bigwedge\Gamma\to\ps\Gamma$
\item[$\mathsf 4$] $\ps\ps \Gamma\to\ps \Gamma$
\item[$\mathsf{Fix}_\ps$] $\ps\Gamma\to\bigwedge_{\gamma\in\Gamma}\ps(\gamma\wedge\ps\Gamma)$
\item[$\mathsf{Ind}_\ps$] Induction for $\ps$: \[\nc\Big(p\to\bigwedge_{\gamma\in\Gamma}\ps(p\wedge \gamma)\Big)\to(p\to\ps\Gamma)\]
\end{description}
\item[]Temporal axioms:
\begin{description}
\item[$\mathsf{Neg}_\next$]$\neg\next p\leftrightarrow\next\neg p$
\item[$\mathsf{And}_\next$] $\next(p\wedge q)\leftrightarrow \next p\wedge \next q$
\item[$\mathsf{Fix}_\hf$] $\hf p\to p\wedge\next\hf p$
\item[$\mathsf{Ind}_\hf$] $\hf(p\to\next p)\to (p\to \hf p)$
\end{description}
\item[$\mathsf{TCont}$]$\ps\next\Gamma\to\next\ps\Gamma$.
\item[]Rules:
\begin{description}
\item[$\mathsf{MP}$] Modus ponens
\item[$\mathsf{Subs}$] Substitution
\item[$\mathsf N_\nc$] Necessitation for $\nc$
\item[$\mathsf N_\next$] Necessitation for $\next$
\item[$\mathsf N_\hf$] Necessitation for $\hf$
\end{description}
\end{description}

\begin{prop}
The above axiomatization is sound for the class of dynamic topological systems.
\end{prop}

\proof
Each axiom and rule above has appeared and shown to be sound either in \cite{kmints} or in \cite{dynamictangle}.
\endproof

Any continuous function satisfies the `tangled' continuity axiom, but it is not logically derivable from the weaker
\[\mathsf{Cont}\phantom{aaaaa}\ps\next p\to\next\ps p,\]
which corresponds to the special case where $\Gamma$ is a singleton \cite{dynamictangle}.

In general we will indicate substitution instances of axioms using parentheses; for example,
\[\mathsf{4}_{\ps}(\varphi)=\ps \ps\varphi\to\ps\varphi.\]

Throughout this paper, $\vdash$ denotes derivability in the system described above, and $\mathcal{DTL}^\ast$ its set of theorems. Note that this strays from the common usage where $\mathcal{DTL}$ is defined semantically, but once we have proven completeness the distinction will become unimportant.

\begin{prop}[Short-term completeness]\label{prelcomp}
Any valid formula in $\mathsf L^\ast_{\ps\next}$ is derivable. Further, the logic over this fragment enjoys the finite model property; that is, if $\varphi\in\mathsf L^\ast_{\ps\next}$ is satisfiable over a dtm, there exists a finite dtm $\mathfrak W$ based on a preorder topology such that $\lb\varphi\rb_\mathfrak W\not=\varnothing.$
\end{prop}

\proof
This is proven in \cite{dynamictangle}.
\endproof

\section{Quasimodels}\label{secndq}

In this section we review a series of results from \cite{me2}, where missing proofs may be found; note, however, that there have been changes to notation and terminology.

We will define quasimodels for $\mathcal{DTL}^\ast$, introduced originally as {\em non-de\-ter\-mi\-nistic quasimodels}. The basic idea is that, while $\mathcal{DTL}^\ast$ is not complete for Kripke models, one can reduce satisfiability in arbitrary dynamical systems to satisfiability in a birelational Kripke structure with certain syntactic constraints. For the construction it is convenient to assign types to worlds, rather than evaluating formulas directly from the propositional variables.

We will denote the set of subformulas of $\varphi$ by ${\mathrm{sub}}(\varphi)$, and define
$${\mathrm{sub}}_{\pm}(\varphi)={\mathrm{sub}}(\varphi)\cup\neg {\mathrm{sub}}(\varphi).$$
We will treat ${\mathrm{sub}}_{\pm}(\varphi)$ as if it were closed under negation, by implicitly identifying $\psi$ with $\neg\neg\psi$. If $\Phi$ is a set of formulas, $\mathrm{sim}(\Phi)$ denotes $\bigcup_{\varphi\in\Phi}\mathrm{sub}(\varphi)$, and $\mathrm{sub}_\pm(\Phi)$ is defined analogously.

\begin{definition}[type]\label{deftype}
A {\em weak type} is a finite set of formulas $\Phi$ such that
\begin{itemize}
\item whenever $\psi\in\Phi$ it follows that $\neg\psi\not\in\Phi$
\item whenever $\psi\wedge\vartheta\in\Phi$, then both $\psi\in\Phi$ and $\vartheta\in\Phi$
\item whenever $\neg(\psi\wedge\vartheta)\in\Phi$, either $\neg\psi\in\Phi$ or $\neg\vartheta\in\Phi$
\item whenever $\hf\psi\in\Phi$, it follows that $\psi\in\Phi$.
\end{itemize}

For a set of formulas $\Psi$, $\Phi\subseteq\mathrm{sub}_\pm(\Psi)$ is a {\em $\Psi$-type} if it is a weak type and, given $\psi\in\mathrm{sub}(\Psi),$ either $\psi\in\Phi$ or $\neg\psi\in\Phi$.

The set of $\Psi$-types will be denoted by $\mathrm{type}(\Psi)$.
\end{definition}

We adopt the general custom of identifying singletons with the element they contain when this does not lead to confusion, so that, for example, we write $\varphi$-type instead of $\cbra\varphi\cket$-type.

\begin{definition}[typed preorder]\label{frame}
Let $S$ be a set preordered by $\peq$.

A {\em weak typing function} on $S$ is a function $t$ which assigns to each $w\in|S|$ a type $t(w)$ such that
\begin{itemize}
\item whenever $w\in S$ and $\ps\Gamma\in t(w)$, there is $v\peq w$ with the property that $\Gamma\subseteq t([v])$\footnote{In other words, for all $\gamma\in\Gamma$ there is $u\sim v$ with $\gamma\in t(u)$.} and
\item whenever $w\in S$, $\neg\ps\Gamma\in t(w)$ and $v\peq w$, there is $\gamma\in\Gamma$ such that $\neg\gamma\in\bigcap t([v])$.
\end{itemize}

If all types are $\Phi$-types, we say $t$ is a {\em $\Phi$-typing function}.

A {\em weakly typed preorder} is a tuple $\mathfrak A=\<|\mathfrak A|,\peq_\mathfrak A,t_\mathfrak A\>$ consisting of a preordered set equipped with a weak typing function; if $t_\mathfrak A$ is a $\Phi$-typing function, then $\mathfrak A$ is a $\Phi$-typed preorder.
\end{definition}

Thus in a $\Phi$-typed preorder, all subformulas of $\Phi$ are decided; on weakly typed structure, only the formulas appearing have a definite value.

\begin{definition}[sensible relation]\label{compatible}
Let $\varphi$ be a formula in $\mathsf L^\ast$ and $\Phi,\Psi$ be finite sets of formulas. The ordered pair $(\Phi,\Psi)$ is {\em sensible} if
\begin{enumerate}
\item whenever $\next\psi\in\Phi$, $\psi\in \Psi$,
\item whenever $\neg\next\psi\in\Phi$, $\neg\psi\in \Psi$,
\item for every formula $\psi$, $\hf\psi\in\Phi$ implies that $\hf\psi\in \Psi$ and
\item for every formula $\psi$, $\ev\psi\in\Phi$ implies that $\psi\in \Phi$ or $\ev\psi\in \Psi$.
\end{enumerate}

Likewise, a pair $(w,v)$ of worlds in a typed preorder $\mathfrak A$ is sensible if $(t(w),t(v))$ is sensible.

A continuous\footnote{By {\em continuous relation} we understand a binary relation $\chi$ such that, whenever $U$ is open, $\chi^{-1}$ is open as well.} relation $\mapsto\subseteq |\mathfrak A|\times |\mathfrak A|$ is {\em sensible} if it is serial and every pair in $\mapsto$ is sensible.

Further, $\mapsto$ is $\omega$-sensible if for all $\ev\psi\in {\mathrm{sub}}_\pm(\varphi)$, if $\ev\psi\in t(w)$, there exist $v\in |\mathfrak A|$ and $N\geq 0$ such that $\psi\in t(v)\text{ and }w\mapsto^Nv.$
\end{definition}

We will refer to formulas of the form $\ev \psi$ as {\em eventualities}. If $\ev\psi\in t_\mathfrak A(w)$, $w\rra^N v$ and $\psi\in t_\mathfrak A(v)$, we will say $v$ {\em realizes} $\ev\psi$, or $\ev\psi$ is realized in time $N$.


\begin{definition}[Quasimodel]\label{ndqm}
Given a finite set of formulas $\Phi$, a {\em $\Phi$-quasimodel} is a tuple
\[\mathfrak{A}=\<|\mathfrak A|,\peq_\mathfrak A,\mapsto_\mathfrak A,t_{\mathfrak A}\>,\]
where $\<|\mathfrak A|,\peq_\mathfrak A,t_{\mathfrak A}\>$ is a $\Phi$-typed Kripke frame and $\mapsto_\mathfrak A$ is an $\omega$-sensible relation on $|\mathfrak A|$.

$\mathfrak A$ {\em satisfies} $\varphi$ if there exists $w_\ast\in |\mathfrak A|$ such that $\varphi\in t_{\mathfrak A}(w_\ast)$.
\end{definition}

We adopt the general practice of dropping subindices when this does not lead confusion, for example writing $\peq$ instead of $\peq_\mathfrak A$.


\section{Generating dynamic topological models from quasimodels}\label{secgen}

Given a $\Phi$-quasimodel $\mathfrak A$, we can construct a dynamic topological model $\lim\mathfrak A$ satisfying the same subformulas of $\Phi$ as $\mathfrak A$; the points of $|\lim\mathfrak A|$ will not be worlds in $|\mathfrak A|$, but rather infinite $\mapsto$-paths.

\subsection{Realizing sequences}
A {\em path} in $\mathfrak A$ is any sequence $\<w_n\>_{n<N}$, with $N\leq \omega$, such that $w_n\mapsto w_{n+1}$.

The continuity of $\rra$ has a natural generalization for finite paths. The following lemma is proven in \cite{me2}:

\begin{lemm}\label{pathcont}
Let $\mathfrak A$ be a $\Phi$-quasimodel,
$\<w_n\>_{n\leq N}$
a finite path and $v_0$ be such that $v_0\peq w_0.$

Then, there exists a path
$\<v_n\>_{n\leq \omega}$
such that, for $n\leq N$, $v_n\peq w_n$.
\end{lemm}

\proof This follows from an easy induction on $N$ using the cotinuity of $\rra$.\endproof

An infinite path
${\vec w}=\< w_n\>_{n<\omega}$
is {\em realizing} if for all $n<\omega$ and $\ev\psi\in t(w_n)$ there exists $K\geq n$ such that $\psi\in t(w_K)$.

Denote the set of realizing paths by $\overrightarrow{|\mathfrak A|}$. Note that $\overrightarrow{|\mathfrak A|}\subseteq|\mathfrak A|^{\mathbb N}$; if we view $|\mathfrak A|$ as a topological space with the preorder topology, then $|\mathfrak A|^{\mathbb N}$ naturally acquires the product topology. Consequently, $\overrightarrow{|\mathfrak A|}$ can be seen as a topological space under the corresponding subspace topology; this topology on $\overrightarrow{|\mathfrak A|}$ will be denoted $\mathcal T_\mathfrak A$.

For $\vec w,\vec v\in \overrightarrow{|\mathfrak A|}$ and $N<\omega$, write $\vec v\stackrel N\peq \vec w$ if $v_n\peq w_n$ for all $n<N$. Then define
\[\mathop{\downarrow_N}(\vec w)=\cbra\vec v\in\overrightarrow{|\mathfrak A|}:\vec v\stackrel N\peq\vec w\cket.\]

Sets of the form $\mathop{\downarrow_N}(\vec w)$ form a basis for $\mathcal T_\mathfrak A$ \cite{me2}.

\subsection{Limit models}

We can define dynamics on $\overrightarrow{|\mathfrak A|}$ by the shift operator $\sigma$, given by
\[\sigma\left(\<w_n\>_{n<\omega}\right)=\<w_{n+1}\>_{n<\omega}.\]
This simply removes the first element in the sequence. The function $\sigma$ is continuous with respect to $\mathcal T_\mathfrak A$.

We can also use $t$ to define a valuation: if $p$ is a propositional variable, set
\[\lb p\rb_{\lim\mathfrak A}=\cbra{\vec w}\in \overrightarrow{|\mathfrak A|}:p\in t\left(w_0\right)\cket.\]

We are now ready to assign a dynamic topological model to every $\varphi$-quasimodel:

\begin{definition}[limit model]
Given a $\Phi$-quasimodel $\mathfrak A$, define
\[\lim \mathfrak{A}=\<\overrightarrow{|\mathfrak A|},\mathcal{T}_{\mathfrak A},\sigma,\lb \cdot\rb_{\lim\mathfrak A}\>\]
to be the {\em limit model} of $\mathfrak{A}$.
\end{definition}

The following was proven in \cite{me2} for monadic formulas (i.e., formulas where $\ps$ is applied only to singletons). The version we present here is a mild generalization.

\begin{lemm}\label{sound}
Suppose $\mathfrak A$ is a $\Phi$-quasimodel, ${\vec w}=\cbra w_n\cket_{n\geq 0}\in |\mathfrak A|$ and $\psi\in {\mathrm{sub}}_{\pm}(\Phi)$.
Then,
\[\<\lim \mathfrak A,{\vec w}\>\models \psi\text{ if and only if } \psi\in t(w_0).\]
\end{lemm}
\proof The proof goes by standard induction of formulas. The induction steps for Boolean operators are trivial, and the steps for the modal operators $\next,\hf$ are covered in \cite{me2}. Hence we will only consider formulas of the form $\ps\Gamma$.

Assume that $\ps\Gamma\in t(w_0)$. In order to prove that $\vec w$ satisfies $\ps\Gamma$, it suffices to show that $\lb\Gamma\rb_{\lim\mathfrak A}$
is tangled in
\[E=\cbra \vec v:\ps\Gamma\in t(v_0)\cket.\]

Let $\vec v\in E$. Suppose $\gamma\in \Gamma$ and $N<\omega$. There is $u_0\peq v_0$ such that $\gamma,\ps\Gamma\in t(u_0)$ (by Definition \ref{frame}.1). But, by Lemma \ref{pathcont}, $u_0$ can then be extended to a realizing path $\vec u\in\mathop\downarrow_N(\vec v)$. Clearly $\vec u\in E$, and by induction on formulas $\vec u\in\lb\gamma\rb_{\lim\mathfrak A}$. Thus every neighborhood of $\vec v$ contains a point in $\lb\gamma\rb_{\lim \mathfrak A}\cap E$. We conclude $E\subseteq \lb\ps\Gamma\rb_{\lim\mathfrak A}$, and since $\vec w\in E$,
\[\<\lim\mathfrak A,\vec w\>\models\ps\Gamma.\]

Now assume $\ps\Gamma\not\in t(w_0)$. 

We must show that there can be no set $E$ containing $\vec w$ such that $\lb\Gamma\rb_{\lim\mathfrak A}$ is tangled in $E$. Here we will use a second induction on $\prec$, in the following sense; if $\vec v$ is any path with $v_0\prec w_0$, we will assume $\ps\Gamma\in t(v_0)$ if and only if $\vec v\in\lb\ps\Gamma\rb_{\lim\mathfrak A}$. Note that, since $\ps\Gamma\not\in t(w_0)$, there is no $v_0\peq w_0$ such that $\ps\Gamma\in t(v_0)$.

Towards a contradiction, suppose that $\lb\Gamma\rb_{\lim\mathfrak A}$ were tangled in $E$. By our induction hypothesis,
\[E\cap\cbra \vec v:v_0\prec w_0\cket=\varnothing.\]

Thus, for all $\vec v\in E\cap\mathop\downarrow_0\vec w$ we have that $v_0\sim w_0$.

But for each $\gamma\in \Gamma$ there must be a point
\[\vec v^\gamma\in E\cap\mathop{\downarrow_0}\vec w\cap \lb\gamma\rb_{\lim\mathfrak A};\] by induction hypothesis (on formulas) this implies that $\gamma\in t(v^\gamma_0)$ and, by the above considerations, $v^\gamma_0\sim w_0$. Thus for each $\gamma\in\Gamma$, $[w_0]$ contains a point $v^\gamma_0$ with $\gamma\in t(v^\gamma_0)$, which by Definition \ref{frame}.2 implies that $\neg\ps\Gamma\not\in t(w_0)$ and thus $\ps\Gamma\in t(w_0)$.
\endproof

This leads us to the main theorem of this section:

\begin{theorem}\label{second}
Let $\varphi$ be a formula of $\mathcal{L}$, and suppose $\varphi$ is satisfied in a quasimodel $\mathfrak A$.

Then, there exists ${{{\vec w}^*}\in |\overrightarrow{\mathfrak A}|}$ such that
\[\<\lim\mathfrak{A},{{\vec w}}^*\>\models\varphi.\]
\end{theorem}
\proof
Pick $w^*\in |\mathfrak A|$ such that $\varphi\in t(w^*)$; $w^*$ can be extended to a realizing path ${\vec w}^*$ \cite{me2}. It follows from Lemma \ref{sound} that
\[\<\lim\mathfrak{A},{\vec w}^*\>\models\varphi.\qedhere\]
\endproof

\section{Simulations}\label{secsim}

We will say that a relation between topological spaces is {\em continuous} if the preimage of any open set is open. Note that this is not the standard definition of continuous relations, which is more involved.

\begin{definition}[Topological simulation]
Let $\mathfrak X$ and $\mathfrak Y$ be topological models. A {\em simulation} is a continuous binary relation
\[\chi\subseteq |\mathfrak X|\times|\mathfrak Y|\]
such that for every $p\in \mathsf{PV}$ and $x\mathrel\chi y$, $x\in\lb p\rb_\mathfrak X$ if and only if $y\in\lb p\rb_\mathfrak Y$.

Given topological models $\mathfrak X$ and $\mathfrak Y$, a point $x\in|\mathfrak X|$ {\em simulates} $y\in|\mathfrak Y|$ if there exists a simulation $\chi\subseteq|\mathfrak X|\times|\mathfrak Y|$ such that $x\mathrel\chi y$; we will write $\<\mathfrak X,x\>\unlhd \<\mathfrak Y,y\>$.
\end{definition}

Note that in the above definition, either $\mathfrak X$ or $\mathfrak Y$ could well be Kripke models, as long as we tacitly identify them with their corresponding preorder topology model.

In the case that both structures are Kripke models then {\em continuity} is the usual `forth' condition for simulations,  namely that, if $v\peq_\mathfrak X w$ and $w\mathrel\chi x$, there is $y\peq_\mathfrak Y x$ such that $v\mathrel \chi y$.

We are also interested in simulations involving typed structures; either between two typed preorders or between a typed preorder and a dtm. We define these below:

\begin{definition}[Typed simulations]
Given two typed preorders $\mathfrak W,\mathfrak V$, we say a {\em simulation} between $\mathfrak W$ and $\mathfrak V$ is a continuous relation $\chi$ such that $w\mathrel \chi v$ implies that $t_\mathfrak W(w)=t_\mathfrak V(v)$.

If $\mathfrak W$ is a typed preorder and $\mathfrak X$ a dtm, $\chi\subseteq|\mathfrak W|\times|\mathfrak X|$ is a {\em simulation} if it is continuous and, whenever $w\mathrel\chi x$, it follows that
\[\<\mathfrak X,x\>\models\bigwedge t_\mathfrak W(w).\]
\end{definition}

\section{The universal state space}\label{secfour}

In this section we define a structure $\mathfrak I(\Phi)$ which we will use to link the semantic framework developed so far with the syntactic constructions we need for a completeness proof. The `worlds' of $\mathfrak I(\Phi)$ are called $\Phi$-states, as defined below. In the end we will be mainly interested in structures $\mathfrak I(\varphi)$ (i.e., when $\Phi$ is a singleton), but it will be convenient to give a more general treatment.

We refer the reader once again to \cite{me2} for omitted proofs.

Roughly, $\Phi$-states are local descriptions of $\Phi$-quasimodels. We will define $\mathrm{len}(\varphi)$ (the {\em lenght} of $\varphi$) as the number of subformulas of $\varphi$, and similarly define $\mathrm{len}(\Phi)$ as $\#\mathrm{sub}(\Phi)$.

\begin{definition}[$\Phi$-state]
Let $\Phi$ be a finite set of formulas. A {\em $\Phi$-state} is a tuple
\[{\mathfrak w}=\<|\mathfrak w|,\peq_\mathfrak w,t_{\mathfrak w},0_{\mathfrak w}\>,\]
consisting of a finite $\Phi$-typed preorder equipped with a distinguished point $0_{\mathfrak w}\in| \mathfrak w|$ satisfying $v\peq_\mathfrak w 0_{\mathfrak w}$ for all $v\in| \mathfrak w|$.
\end{definition}

In other words, a $\Phi$-state is a $\Phi$-typed, finite, pointed Kripke frame, or a local Kripke frame from \cite{me2}. We will write $t( \mathfrak w)$ instead of $t_{\mathfrak w}\left(0_{\mathfrak w}\right)$. As always, we will include the subindex in $\peq_\mathfrak w$ only when necessary. {\em Weak states} are defined analogously but based on weakly typed preorders.

\begin{definition}[norm of a state]
Given a state $ \mathfrak w$, we define $\mathrm{hgt}( \mathfrak w)$ as the maximum length of a sequence of worlds
\[w_0\prec w_1\prec w...\prec w_N\]
such that $w_n\in| \mathfrak w|$ for all $n\leq N$. Similarly, $\mathrm{wdt}( \mathfrak w)$ is defined as the maximum $N$ such that there exist $w\in| \mathfrak w|$ with $N$ $\peq$-incomparable daughters.

We then define the {\em norm} of $ \mathfrak w$, denoted $\mathrm{nrm}\left( \mathfrak w\right)$, by
\[\mathrm{nrm}\left( \mathfrak w\right)=\max(\mathrm{hgt}( \mathfrak w),\mathrm{wdt}( \mathfrak w)).\]
\end{definition}

We will make the assumption that no two worlds of $ \mathfrak w$ are indistinguishable; that is, if $w\sim v$ and $t_{ \mathfrak w}(w)=t_{ \mathfrak w}(v)$, then $w=v$. This will immediately bound the size of each cluster by $2^{\mathrm{len}(\Phi)}$. Thus bounding the height and width of $ \mathfrak w$ gives us a bound on $\#| \mathfrak w|$ (and vice-versa). In particular, it follows that there are only finitely many $\Phi$-states with a given norm.

Many times it will be useful to compare different $\Phi$-states and express relations between them. These relations will appear throughout the rest of the text.

\begin{definition}
Say that ${\mathfrak w}$ {\em simulates} ${\mathfrak v}$, denoted ${\mathfrak w}\unlhd{\mathfrak v}$, if there exists a simulation $\chi\subseteq|\mathfrak w|\times |\mathfrak v|$ such that $0_{\mathfrak w}\mathrel\chi 0_{\mathfrak v}$.
\end{definition}

\begin{definition}[substate]
Say that $\mathfrak v$ is a {\em substate} of $\mathfrak w$, written $\mathfrak v\peq\mathfrak w$, if $0_\mathfrak v\in|\mathfrak w|$ and $\mathfrak v$ is a generated substructure of $\mathfrak w$.
\end{definition} 

Below, $\mathrm{sub}_\Diamond(\Psi)$ denotes the set of all formulas $\ps\Gamma\in\mathrm{sub}(\Psi)$.

\begin{definition}[termporal successor]\label{ts}
Say $ \mathfrak w$ is a {\em temporal successor} of $ \mathfrak v$, denoted $ \mathfrak w\mapsto \mathfrak v$, if there exists a sensible relation
$g\subseteq |\mathfrak w|\times |\mathfrak v|$
such that $0_{\mathfrak w}\mathrel g0_{\mathfrak v}$.

If $\mathfrak w\mapsto\mathfrak v$ and
\[\mathrm{nrm}(\mathfrak v)\leq \mathrm{nrm}(\mathfrak w)+\#\bigcup_{w\in|\mathfrak w|}\mathrm{sub}_\Diamond(t_\mathfrak w(w)),\]
we will write $\mathfrak w\dot\mapsto\mathfrak v$ and say $\mathfrak v$ is a {\em small} temporal successor of $\mathfrak w$.
\end{definition}

Note that if $\mathfrak w$ is a $\Phi$-type, $\bigcup_{w\in|\mathfrak w|}\mathrm{sub}_\Diamond(t_\mathfrak w(w))$ becomes $\mathrm{sub}_\ps(\Phi)$.

With these relations in mind, the class of finite $\Phi$-states can be viewed as a typed structure on its own right.

\begin{definition}[$\mathfrak I_K(\Phi)$]
Let $\Phi$ be a finite set of $\mathsf L^\ast$-formulas and $K\geq 0$. Define $|\mathfrak I_K(\Phi)|$ to be the set of $\Phi$-states $ \mathfrak w$ such that
\[\mathrm{nrm}( \mathfrak w)\leq (K+1)\mathrm{len}(\Phi).\]
\end{definition}

Now, consider
\[|\mathfrak I(\Phi)|=\bigcup_{k<\omega}|\mathfrak I_k(\Phi)|\]
(evidently this is the set of all finite $\Phi$-states).
Define
\[\mathfrak{I}(\Phi)=\<|\mathfrak I(\Phi)|,\peq,\rra,t\>,\]
where $t( \mathfrak w)=t_{ \mathfrak w}(0_{\mathfrak w})$.

We then have the following restatement of a result in \cite{me2}:

\begin{prop}\label{total}
Let $\Phi$ be a finite set of formulas. Then,
\begin{enumerate}
\item If $\mathfrak w$ is any $\Phi$-state, there exists $\mathfrak v\in|\mathfrak I_0(\Phi)|$ such that $\mathfrak v\unlhd\mathfrak w$ and
\item If $\mathfrak w$ is a weak state and $\mathfrak w\mapsto\mathfrak v\in|\mathfrak I(\Phi)|$, there is $\mathfrak u\unlhd\mathfrak v$ such that $\mathfrak w\dot\mapsto\mathfrak u$. 
\end{enumerate}
\end{prop}

\proof
We omit the proof. It proceeds by induction on the height of a state, deleting worlds until we reach a model of the desired size.
\endproof

While $\mathfrak I(\Phi)$ is a $\Phi$-typed preorder with a sensible relation $\rra$, it is not necessarily $\omega$-sensible, so $\mathfrak{I}(\Phi)$ is not a quasimodel. Nevertheless, it will be very useful as a universal structure. In particular, if $\varphi$ is satisfiable, it can be satisfied in a quasimodel which is a substructure of $\mathfrak I(\varphi)$.

More specifically, define a non-empty set $U\subseteq|\mathfrak I(\Phi)|$ to be {\em regular} if $U$ is open (i.e., downward-closed under $\peq$) and $\mapsto\upharpoonright U$\footnote{We denote restriction by $\upharpoonright$, so that for example $\mapsto\upharpoonright U=\mapsto\cap(U\times U)$.} is $\omega$-sensible. The following should then be fairly obvious from Definition \ref{ndqm}:

\begin{lemm}\label{regquasi}
If $U$ is a regular subset of $|\mathfrak I(\Phi)|$, then $\mathfrak I(\Phi)\upharpoonright U$ is a $\Phi$-quasimodel.
\end{lemm}

We will refer to quasimodels of the form $\mathfrak I(\Phi)\upharpoonright U$, with $U$ regular, as {\em regular quasimodels}. Much of what follows will be devoted to defining a `canonical quasimodel' for a given consistent set of formulas $\Phi$, and this quasimodel will be regular.

\section{Simulation formulas}\label{secsubchar}

The primary motivation for extending the language of $\mathcal{DTL}$ to use a polyadic modality is that, for our completeness proof, it is essential to be able to define simulability by finite $\Phi$-states. This cannot be done in the standard modal language, but in the extended language the situation is different \cite{me:simulability}. In this section, we will discuss the formulas $\mathrm{Sim}(\mathfrak w)$, which define the property being simulated by $\mathfrak w$.

We generalize the notion of substitution to states as follows: if $\mathfrak w$ is a state, $\vec p$ a tuple of variables and $\vec\psi$ a tuple of formulas, we write $\mathfrak w[\vec p/\vec\psi]$ as the state $\mathfrak v$ which is identical to $\mathfrak w$ except that, for $w\in|\mathfrak w|$, we have
\[t_\mathfrak v(w)=\cbra \delta[\vec p/\vec\psi]:\delta\in t_\mathfrak w(w)\cket.\]



For a state $\mathfrak w$, define $\mathfrak w^{\mathsf p}$ as the state which is identical to $\mathfrak w$ except that $t_{\mathfrak w^{\mathsf p}}$ is given as follows: for each type $\Phi$ in the range of $t_{\mathfrak w}$ (which we denote $\mathrm{rng}(t_\mathfrak w)$), we introduce a new propositional variable $p_\Phi$.

Then we put
\[t_{\mathfrak w^{\mathsf p}}(w)=\cbra p_{t_\mathfrak w(w)}\cket\cup\cbra \neg p_\Psi:\Psi\not=t_\mathfrak w(w)\cket.\]

Say a weak state $\mathfrak w$ is {\em distinctly typed} if, whenever $t_\mathfrak w(w)\not=t_\mathfrak w(v)$, there is $\psi\in t_\mathfrak w(w)$ such that $\neg\psi\in t_\mathfrak w(v)$ (or vice-versa). Note that $\Phi$-typed states are distinctly typed.

We then get:

\begin{prop}\label{simulability}
Given a distinctly typed state $\mathfrak w$, there exists a formula $\mathrm{Sim}(\mathfrak w)$ such that, for every topological model $\mathfrak X$ and $x\in|\mathfrak X|$, $x\in\lb\mathrm{Sim}(\mathfrak w)\rb_\mathfrak X$ if and only if $\mathfrak w\unlhd\<\mathfrak X,x\>$.

Further, we can define $\mathrm{Sim}$ so that
\[\mathrm{Sim}(\mathfrak w)=\mathrm{Sim}(\mathfrak w^{\mathsf p})\bra p_\Psi/\bigwedge\Psi\ket_{\Psi\in\mathrm{rng}(t_\mathfrak w)}.\]
\end{prop}

\proof
Let $\mathfrak w$ be a distinctly typed state.

In \cite{me:simulability} it is shown that there exists a formula $\mathrm
{Sim}(\mathfrak w^\mathsf p)$ such that, for any topological model $\mathfrak X$ and $x\in|\mathfrak X|$, $\mathfrak w^{\mathsf p}\unlhd \<\mathfrak X,x\>$ if and only if $\<\mathfrak X,x\>\models\vartheta$. We then set
\[\mathrm{Sim}(\mathfrak w)=\mathrm{Sim}(\mathfrak w^{\mathsf p})\bra p_\Psi/\bigwedge\Psi\ket_{\Psi\in\mathrm{rng}(t_\mathfrak w)},\]
as indicated.

Now, if $\mathfrak Y$ is any dtm, let $\mathfrak Y^{\mathsf p}$ be the dtm obtained by extending $\lb\cdot\rb_\mathfrak Y$ so that $\lb p_\Psi\rb_{\mathfrak Y^{\mathsf p}}=\lb \bigwedge\Psi\rb_\mathfrak Y$. Note that these sets are disjoint because different types are mutually inconsistent, given that $\mathfrak w$ is distinctly typed. 

We then have that
\[\begin{array}{lcl}
\<\mathfrak Y,x\>\models\mathrm{Sim}(\mathfrak w)&\Leftrightarrow&\<\mathfrak Y^{\mathsf p},x\>\models\mathrm{Sim}(\mathfrak w^{\mathsf p})\\\\
&\Leftrightarrow&\mathfrak w^{\mathsf p}\unlhd\<\mathfrak Y^{\mathsf p},x\>\\\\
&\Leftrightarrow&\mathfrak w\unlhd\<\mathfrak Y,x\>,
\end{array}\]
as claimed.
\endproof



\begin{prop}\label{propsub}
Simulation formulas have the following properties:
\begin{enumerate}
\item If $\psi\in t(\mathfrak w)$, then $\vdash\mathrm{Sim}(\mathfrak w)\to\psi$;
\item if $\mathfrak v\unlhd\mathfrak w$ then $\vdash \mathrm{Sim}(\mathfrak w)\to\mathrm{Sim}(\mathfrak v)$;
\item if $\mathfrak v\peq\mathfrak w$ then $\vdash \mathrm{Sim}(\mathfrak w)\to\ps\mathrm{Sim}(\mathfrak v)$;
\item if $\psi\in\mathrm{sub}(\Phi)\cup\neg\mathrm{sub}(\Phi)$,
\[\vdash\psi\to\bigvee\cbra \mathrm{Sim}(\mathfrak w):\mathfrak w\in\mathfrak I_0(\Phi)\text{ and }\psi\in t(\mathfrak w)\cket;\]
\item for all $\mathfrak w\in\mathfrak I(\Phi)$,
\[\vdash\mathrm{Sim}(\mathfrak w)\to\next\bigvee_{\mathfrak w\stackrel\cdot\mapsto\mathfrak v}\mathrm{Sim}(\mathfrak v).\]
\end{enumerate}
\end{prop}

\proof
As we shall see, these claims are consequences of Propositions \ref{prelcomp}, \ref{total} and \ref{simulability}. What we shall do is show that each of these formulas is a (subsitution instance of an) $\mathsf{S4C}^\ast$ validity, which by Proposition \ref{prelcomp} implies that it is derivable.

\paragraph{1.} Suppose $\psi\in t(\mathfrak w)$. Any point on a dtm satisfying $\mathrm{Sim}(\mathfrak w)$ must satisfy $t(\mathfrak w)$ (by the definition of a simulation), which in this case includes $\psi$. Hence $\mathrm{Sim}(\mathfrak w)\to\psi$ is a substitution instance of an $\mathsf{S4C}^\ast$ validity, as claimed.

\paragraph{2.} This expresses the transitivity of $\unlhd$. Namely, suppose $\<\mathfrak X,x\>\models \mathrm{Sim}(\mathfrak w)$ (so that $\mathfrak w\unlhd\<\mathfrak X,x\>$) and $\mathfrak v\unlhd\mathfrak w$. Then, clearly $\mathfrak v\unlhd\<\mathfrak X,x\>$ (by composing the itermediate simulations) so by Proposition \ref{simulability} we have that $\<\mathfrak X,x\>\models\mathrm{Sim}(\mathfrak v)$.

\paragraph{3.} Suppose $\mathfrak v\peq\mathfrak w$ and $\<\mathfrak X,x\>$ is a dtm satisfying $\mathrm{Sim}(\mathfrak w)$, so that there is a simulation $\chi$ between $\mathfrak w$ and $\<\mathfrak X,x\>$ with $0_\mathfrak w\mathrel\chi x$. Since simulations are continuous and the only neighborhood of $0_\mathfrak w\in|\mathfrak w|$ is all of $|\mathfrak w|$, it follows that, given a neighborhood $U$ of $x$, $\chi^{-1}(U)=|\mathfrak w|$.

In particular, $0_\mathfrak v\in\chi^{-1}(U)$, so that $\zeta=\chi\upharpoonright|\mathfrak v|$ is a simulation between $\mathfrak v$ and $\mathfrak X$ with $0_\mathfrak v\mathrel \zeta y$ for some $y\in U$. Thus by Proposition \ref{simulability}, $y$ satisfies $\mathrm{Sim}(\mathfrak v)$ and, since $U$ was arbitrary, $x\in\overline{\lb\mathrm{Sim}(\mathfrak v)\rb}_\mathfrak X$ or, equivalently, $x$ satisfies $\ps\mathrm{Sim}(\mathfrak v)$.

\paragraph{4.}
It will be useful to define a variant of $\mathfrak w$, which we denote $\mathfrak w^{\mathsf q}$. For this, add a new propositional variable $q_\delta$ for each $\hf\delta\in\mathrm{sub}(\Phi)$. Given a formula $\gamma$, let $\gamma^\mathsf q$ be the result of replacing each outermost occurrence of $\hf\delta$ in $\gamma$ by $q_\delta$. Similarly, if $\Gamma$ is a set of formulas, $\Gamma^\mathsf q$ denotes the set $\cbra \gamma^\mathsf q:\gamma\in \Gamma\cket$ and $\mathfrak w^{\mathsf q}$ is the state which is identical to $\mathfrak w$ except that $t_{\mathfrak w}(w)$ is replaced by $t_\mathfrak w^\mathsf q(w)$.

We claim that
\[\mathsf{S4C}^\ast\vdash\eta,\]
where $\eta$ is the formula
\[\psi^\mathsf q\to \bigvee\cbra \mathrm{Sim}(\mathfrak w):\mathfrak w\in\mathfrak I_0(\mathrm {sub}^\mathsf q(\Phi))\text{ and }\psi^\mathsf q\in t(\mathfrak w)\cket.\]
By the finite model property for $\mathsf{S4C}^\ast$ (Propostion \ref{prelcomp}), it suffices to check that this formula is valid on every finite dynamic Kripke model $\mathfrak W$.

Suppose, then, that $\mathfrak W$ is finite and $\<\mathfrak W,w\>\models\psi^\mathsf q$. We obtain a $\mathrm{sub}^\mathsf q(\Phi)$-state\footnote{Note that $\mathrm{sub}^\mathsf q(\Phi)$ may not be equal to $\mathrm{sub}(\Phi^\mathsf q)$! Consider, for example, $\Phi=\cbra\hf p\cket$.} $\mathfrak v$ by letting $|\mathfrak v|=\mathop\downarrow w$, $0_\mathfrak v=w$ and
\[t_\mathfrak v(v)=\cbra\vartheta\in \mathrm{sub}_\pm^\mathsf q(\Phi):v\in\lb\vartheta\rb_\mathfrak W\cket.\]

Now, by Proposition \ref{total} there is $\mathfrak w\in\mathfrak I_0(\mathrm{sub}^\mathsf q(\Phi))$ such that $\mathfrak w\unlhd\mathfrak v$ and thus $w\in\lb\mathrm{Sim}(\mathfrak w)\rb_{\mathfrak X^\mathsf q}$. Hence $\<\mathfrak W,w\>$ satisfies
\[\bigvee\cbra \mathrm{Sim}(\mathfrak w):\mathfrak w\in\mathfrak I_0(\mathrm{sub}^\mathsf q(\Phi))\text{ and }\psi^\mathsf q\in t(\mathfrak w)\cket,\]
and $\eta^\mathsf q$ is valid.

Let us write $[\vec q/\hf\vec\delta]$ instead of $\Big[ q_\delta/\hf\delta\Big]_{\hf\delta\in\mathrm{sub}\Phi}$. We then have that $\eta[\vec q/\hf\vec\delta]$ is equal to
\begin{equation}\label{etaquation}
\psi\to\bigvee\cbra \mathrm{Sim}(\mathfrak w)[\vec q/\hf\vec\delta]:\mathfrak w\in\mathfrak I_0(\mathrm {sub}^\mathsf q(\Phi))\text{ and }\psi^\mathsf q\in t(\mathfrak w)\cket;
\end{equation}
this is derivable by applying $\mathsf{Subs}$ to $\eta$, and very close to our goal, except that it may be that for some $\mathrm {sub}^\mathsf q(\Phi)$-state $\mathfrak w$ with $\psi^\mathsf q\in t(\mathfrak w)$, the structure $\mathfrak w[\vec q/\hf\vec\delta]$ is not a $\Phi$-state.

In other words, it may be the case that for some $v\in|\mathfrak w|$, $t_{\mathfrak w}(v)[\vec q/\hf\vec\delta]$ is not a $\Phi$-type. However, we shall see that in such cases, $\mathrm{Sim}(\mathfrak w)[\vec q/\hf\vec\delta]$ is inconsistent, and we can remove it from the disjunction.

Looking at Definition \ref{deftype}, the only property that may fail is that $\hf\delta\in t_\mathfrak w(v)$ but $\delta\not\in t_\mathfrak w(v)$. In this case, it follows that $\neg\delta\in t_\mathfrak w(v)$; this is because $t_\mathfrak w(v)$ is a $\mathrm{sub}^\mathsf q(\Phi)$-type and $\delta^\mathsf q\in \mathrm{sub}^\mathsf q(\Phi)$, so given that $\delta^\mathsf q\not\in t_\mathfrak w(v)$, we have $\neg\delta^\mathsf q\in t_\mathfrak w(v)$ and thus $\neg\delta\in t_\mathfrak w(v)[\vec q/\hf\vec\delta].$

On the other hand, $\vdash\hf\delta\to\delta$, so
\[\vdash \bigwedge t_\mathfrak w(v)\to(\delta\wedge\neg\delta),\]
i.e. $t_\mathfrak w(v)$ is inconsistent.

Now, let $\mathfrak v\peq\mathfrak w$ be the substructure with $0_\mathfrak v=v$. We have, by Item 3, that
\[\vdash\mathrm{Sim}(\mathfrak w)\to\ps\mathrm{Sim}(\mathfrak v),\]
while by Item 1,
\[\vdash\mathrm{Sim}(\mathfrak v)\to\bigwedge t_\mathfrak w(v);\]
together, these imply that $\vdash\neg\mathrm{Sim}(\mathfrak w).$

Thus we can remove from (\ref{etaquation}) all disjoints of the form $\mathrm{Sim}(\mathfrak w)[\vec q/\hf\vec\delta]$ where $\mathfrak w[\vec q/\hf\vec\delta]\not\in |\mathfrak I_0(\Phi)|$ and obtain the desired result.

\paragraph{5.}

Given a type $\Psi$, define $\Psi^+$ by
\[\Psi^+=\Psi\cup\cbra \next\hf\psi:\hf\psi\in \Phi\cket\cup \cbra\next\ev\psi:\ev\psi\in\Psi\text{ but }\psi\not\in\Psi\cket.\]

It should be fairly clear that
\begin{equation}\label{henceeq}
\mathcal{DTL}^\ast\vdash \bigwedge\Psi\leftrightarrow\bigwedge\Psi^+.
\end{equation}

Define $\mathfrak w^+$ to be the structure which is identical to $\mathfrak w$ except that $t_{\mathfrak w^+}=t_{\mathfrak w}^+$. Note that $\mathfrak w^+$ is, in general, not a proper $\Phi$-state because we have added some formulas which may not be subformulas of $\Phi$; however, it is distincly typed.

Thus we can use Proposition \ref{simulability} to see that
\[\mathrm{Sim}(\mathfrak w)=\mathrm{Sim}(\mathfrak w^{\mathsf p})\big[p_\Psi/\bigwedge\Psi\Big ]_{\Psi\in\mathrm{rng}(t_\mathfrak w)},\]
\[\mathrm{Sim}(\mathfrak w^+)=\mathrm{Sim}(\mathfrak w^{\mathsf p})\Big[p_\Psi/\bigwedge\Psi^+\Big ]_{\Psi\in\mathrm{rng}(t_\mathfrak w)}\]
and by (\ref{henceeq}),
\begin{equation}\label{plusto}
\mathcal{DTL}^\ast\vdash \mathrm{Sim}(\mathfrak w)\leftrightarrow\mathrm{Sim}(\mathfrak w^+).\footnote{We are also using the rule
\[\dfrac{\alpha\leftrightarrow\beta}{\psi[p/\alpha]\leftrightarrow\psi[p/\beta]},\]
which is easily shown to be admissible.
}
\end{equation}

As before, we use the finite model property of $\mathsf{S4C}^\ast$. Suppose that $\mathfrak W$ is a finite dtm and $\<\mathfrak W,w\>\models \mathrm{Sim}(\mathfrak w^{+\mathsf q})$ (where we are writing $\mathfrak w^{+\mathsf q}$ instead of $(\mathfrak w^+)^\mathsf q$, defined in the previous item). As before, we can see $\mathfrak W\upharpoonright \mathop\downarrow f_\mathfrak W(w)$ as a $\mathrm{sub}^\mathsf q(\Phi)$-state, which we shall call $\mathfrak u$.

Then, by Proposition \ref{simulability}, $\mathfrak w^{+\mathsf q}\unlhd \<\mathfrak W,w\>$, and composing with $f_\mathfrak W$ we obtain that $\mathfrak w^{+\mathsf q}\mapsto \mathfrak u$. By Proposition \ref{total} there exists $\mathfrak v$ such that $\mathfrak v\unlhd\mathfrak u$ and $\mathfrak w^{+\mathsf q}\dot\mapsto\mathfrak v$.

It follows that $w$ satisfies
\[\next\bigvee_{\mathfrak w^{+\mathsf q}\dot\mapsto\mathfrak v}\mathrm{Sim}(\mathfrak v);\]
since $\<\mathfrak W,w\>$ was arbitrary,
\[\mathsf{S4C}^\ast\vdash\mathrm{Sim}(\mathfrak w^{+\mathsf q})\to\next\bigvee_{\mathfrak w^{+\mathsf q}\dot\mapsto\mathfrak v}\mathrm{Sim}(\mathfrak v),\]
given that the latter is a formula in $\mathsf L^\ast_{\ps\next}$ and $\mathsf{S4C}^\ast$ is complete.

But then, using $\mathsf{Subs}$, it follows that
\[
\mathcal{DTL}^\ast\vdash \Big(\mathrm{Sim}(\mathfrak w^{+\mathsf q})\to\next\bigvee_{\mathfrak w^{+\mathsf q}\dot\mapsto\mathfrak v}\mathrm{Sim}(\mathfrak v)\Big)[\vec q/\hf\vec\delta],\]
or equivalently
\begin{equation}\label{eqhere}
\mathcal{DTL}^\ast\vdash\mathrm{Sim}(\mathfrak w^{+})\to\next\bigvee_{\mathfrak w^{+\mathsf q}\dot\mapsto\mathfrak v}\mathrm{Sim}\left(\mathfrak v[\vec q/\hf\vec\delta]\right).
\end{equation}

We claim now that if $\mathfrak w^{+\mathsf q}\dot\mapsto\mathfrak v$, then
\[\mathfrak w\dot\mapsto\mathfrak v[\vec q/\hf\vec\delta].\]

Note that for any set of formulas $\Psi$,
\[\#\mathrm{sub}_\ps(\Psi)=\#\mathrm{sub}_\ps(\Psi^+)=\#\mathrm{sub^\mathsf q_\ps(\Psi^+)};\]
thus if $\mathfrak w^{+\mathsf q}\dot\mapsto\mathfrak v$ and $\mathfrak w \mapsto\mathfrak v[\vec q/\hf\vec\delta],$ we automatically have that $\mathfrak v[\vec q/\hf\vec\delta]$ is a {\em small} temporal successor of $\mathfrak w$. Hence we need only prove that $\mathfrak v[\vec q/\hf\vec\delta]$ is a temporal successor of $\mathfrak w$.

In fact, if we have a sensible relation $g$ between $\mathfrak w^{+\mathsf q}$ and $\mathfrak v$, $g$ is still sensible between $\mathfrak w$ and $\mathfrak v[\vec q/\hf\vec\delta].$ To see this, assume that $w\mathrel g v$; let us check that the pair $(t_\mathfrak w(w),t_\mathfrak v(v))$ is sensible (see Definition \ref{compatible}).

\begin{enumerate}
\item Suppose that $\next\psi\in t_{\mathfrak w^+}(w)$. Then, $\psi$ is of the form $\gamma[\vec q/\hf\vec\delta],$ so that $\next\gamma\in  t_{\mathfrak w^{+\mathsf q}}(w)$ and thus $\gamma\in t_{\mathfrak w^{+\mathsf q}}(v)$; it follows that $\psi\in t_{\mathfrak v[\vec q/\hf\vec\delta]}(v)$, as required.

\item The case for $\neg\next\psi$ is similar.

\item If $\hf\psi\in t_{\mathfrak w}(w)$ and $w\mathrel g v$, by construction $\next\hf\psi\in t^+_{\mathfrak w}(w)$, so that $\next q_\psi\in t_{\mathfrak w^{+\mathsf q}}(w)$. This implies that $q_\psi\in t_{\mathfrak v}(v)$; the latter in turn implies that $\hf\psi\in t_{\mathfrak v[\vec q/\hf\vec\delta]}(v),$ which is what we wanted. 

\item The condition for $\ev\psi$ is similar and we skip it.
\end{enumerate}
Thus we can replace (\ref{eqhere}) by
\begin{equation}\label{eqthere}
\mathcal{DTL}^\ast\vdash\mathrm{Sim}(\mathfrak w^{+})\to\next\bigvee_{\mathfrak w\dot\mapsto\mathfrak v}\mathrm{Sim}(\mathfrak v).
\end{equation}

Putting together (\ref{plusto}) and (\ref{eqthere}) we get that
\[
\mathcal{DTL}^\ast\vdash\mathrm{Sim}(\mathfrak w)\to\next\bigvee_{\mathfrak w\dot\mapsto\mathfrak v}\mathrm{Sim}(\mathfrak v).\qedhere
\]
\endproof

\section{Canonical quasimodels}\label{seccan}

We are now ready to define our canonical quasimodels. Given a finite set of formulas $\Phi$, we shall define a quasimodel $\mathfrak Q(\Phi)$ satisfying all consistent $\Phi$-types. This quasimodel shall be a substructure of $\mathfrak I(\Phi)$, however, it will only contain states which are {\em consistent} in the following sense:

\begin{definition}[Consistent states]\label{defsound}
We say a state $\mathfrak w$ is {\em inconsistent} if $\vdash\neg\mathrm{Sim}(\mathfrak w)$; otherwise, it is {\em consistent}.

We will denote the set of consistent $\Phi$-states by $\mathsf{Cons}(\Phi)$.
\end{definition}

With this we are ready to define our canonical quasimodels. However, as showing that they are indeed quasimodels will take some work, we shall baptize them, for now, as {\em canonical structures}:

\begin{definition}[Canonical structures]
Given a set of formulas $\Phi$, we define the {\em canonical strucure} for $\Phi$, denoted $\mathfrak Q(\Phi)$, as $\mathfrak I(\Phi)\upharpoonright \mathsf{Cons}(\Phi)$.

More specifically, we define
\[\mathfrak Q(\Phi)=\langle|\mathfrak Q(\Phi)|,\peq_{\mathfrak Q(\Phi)},\mapsto_{\mathfrak Q(\Phi)},t_{\mathfrak Q(\Phi)}\rangle\] by
\[
\begin{array}{lcl}
|\mathfrak Q(\Phi)|&=&\mathsf{Cons}(\Phi)\\
\peq_{\mathfrak Q(\Phi)}&=&\peq_{\mathfrak I(\Phi)}\cap (\mathsf{Cons}(\Phi)\times \mathsf{Cons}(\Phi))\\
\mapsto_{\mathfrak Q(\Phi)}&=&\mapsto_{\mathfrak I(\Phi)}\cap (\mathsf{Cons}(\Phi)\times \mathsf{Cons}(\Phi))\\
t_{\mathfrak Q(\Phi)}&=&t_{\mathfrak I(\Phi)}\cap (\mathsf{Cons}(\Phi)\times 2^{\mathsf L^\ast})\\
\end{array}
\]
\end{definition}

Our strategy from here on will be to show that canonical structures are indeed quasimodels; once we establish this, completeness of $\mathcal {DTL}^\ast$ is an easy consequence. The most involved step will be showing that $\mapsto_{\mathfrak Q(\Phi)}$ is $\omega$-consistent; however, we may already prove some nice properties of $\mathfrak Q(\Phi)$.

\begin{lemm}\label{niceprop}
If $\Phi$ is a finite set of formulas, $|\mathfrak Q(\Phi)|$ is an open subset of $|\mathfrak I(\Phi)|$ and $\mapsto_{\mathfrak Q(\Phi)}$ is serial.
\end{lemm}

\proof

To check that $|\mathfrak Q(\Phi)|$ is open, let $\mathfrak w\in |\mathfrak Q(\varphi)|$ and suppose $\mathfrak v\peq\mathfrak w$. Now, by Proposition \ref{propsub}.3, we have that
\[\vdash\mathrm{Sim}(\mathfrak w)\to\ps\mathrm{Sim}(\mathfrak v);\]
hence if $\mathfrak w$ is consistent, so is $\mathfrak v$.

To see that $\mapsto_{\mathfrak Q(\Phi)}$ is serial, observe that by Proposition \ref{propsub}.5, if  $\mathfrak w\in |\mathfrak Q(\varphi)|$
for all $\mathfrak w\in\mathfrak I(\Phi)$,
\[\vdash\mathrm{Sim}(\mathfrak w)\to\next\bigvee_{\mathfrak w\stackrel\cdot\mapsto\mathfrak v}\mathrm{Sim}(\mathfrak v);\]
since $\mathfrak w$ is consistent, it follows that for some $\mathfrak v$ with $\mathfrak w\dot\mapsto\mathfrak v$, $\mathfrak v$ is consistent as well, and thus $\mathfrak v\in|\mathfrak Q(\Phi)|.$
\endproof

\section{Efficiency}\label{re}
One of the primary difficulties in the study of $\mathcal{DTL}$ is that we must consider an infinite number of states, so that one cannot tell {\em a priori} how long it will take for a formula of the form $\ev\psi$ to be realized. However, we can remedy this using ideas from \cite{pml,konev} which we elaborate below.

If $\mathfrak A$ is a $\Phi$-quasimodel, an {\em eventuality} is any formula of the form $\ev\psi\in {\mathrm{sub}}_\pm(\Phi)$.

\begin{definition}[efficiency]
Let $\vec{\mathfrak w}=\<\mathfrak w_n\>_{n\leq N}$ be a finite path of $\Phi$-states and $\ev\psi\in t(\mathfrak w)$.

An {\em inefficiency} in $\vec{\mathfrak w}$ is a pair
$M_1<M_2$ such that ${\mathfrak {w}}_{M_1}\unlhd{\mathfrak {w}}_{M_2}$.

The path $\vec{\mathfrak w_n}$ is {\em efficient} if
\begin{enumerate}
\item for all $n<N$, $\mathfrak w_n\dot\mapsto\mathfrak w_{n+1}$ and
\item it contains no inefficiencies.
\end{enumerate}
\end{definition}

Roughly, the previous definition says that if a path realizes an eventuality efficiently, no state should be `repeated'. Otherwise, the path between them gives us a sort of loop which we could simply skip.

Efficient paths are very useful because, while there may be infinitely many paths beginning on a $\Phi$-state $\mathfrak w$, there are only finitely many efficient ones.

The following is a restatement of a result from \cite{pml} which was first applied to $\mathcal{DTL}$ in \cite{konev}:

\begin{lemm}\label{finitist}
Given a $\Phi$-state $\mathfrak w$, there are only finitely many efficient paths $\vec{\mathfrak v}$ with $\mathfrak v_0=\mathfrak w$.
\end{lemm}

\proof
This is a consequence of Kruskal's tree theorem together with K\"onig's Lemma; the proof may be found in \cite{me2,pml,konev}.
\endproof

\section{$\omega$-Sensibility}\label{tss}

In this section we shall show that $\mapsto_{\mathfrak Q(\Phi)}$ is $\omega$-sensible, the most difficult step in proving that $\mathfrak Q(\Phi)$ is a quasimodel. In other words, we must show that, given $\mathfrak w\in|\mathfrak Q(\Phi)|$ and $\ev\psi\in t(\mathfrak w)$, there is a finite path
\[\mathfrak w=\mathfrak w_0\mapsto\mathfrak w_1\mapsto\hdots\mathfrak w_N,\]
where $\psi\in t(\mathfrak w_N)$ and each $\mathfrak w_n\in|\mathfrak Q(\Phi)|$.

Another way of saying this is that there should be some $\mathfrak v\in|\mathfrak Q(\Phi)|$ which is {\em reachable} from $\mathfrak w$ with $\psi\in t(\mathfrak v)$.

Here is where our notion of efficiency will be useful; for, in fact, we need only focus our attention on those $\mathfrak v$ which are reachable from $\mathfrak w$ via an {\em efficient} path.

\begin{definition}[Reachability]
Let $\mathfrak w$ be a $\Phi$-state.

Say a $\Phi$-state $\mathfrak v$ is {\em reachable} from $\mathfrak w$ if there is a finite efficient path 
\[\vec{\mathfrak u}=\<\mathfrak u_0,...,\mathfrak u_N\>\]
of consistent states with $\mathfrak u_0=\mathfrak w$ and $\mathfrak u_N=\mathfrak v$.

Let $\rho(\mathfrak w)$ be the set of all states that are reachable from $\mathfrak w$.
\end{definition}

This notion of reachability is very convenient because of the following:

\begin{lemm}\label{finite2}
Given $\mathfrak w\in|\mathfrak Q(\Phi)|$, $\rho(\mathfrak w)$ is finite.
\end{lemm}

\proof
This is an immediate consequence of Lemma \ref{finitist}.
\endproof

Now, we shall see that to check that $\mapsto_{\mathfrak Q(\Phi)}$ is $\omega$-sensible, we need only consider efficient paths. But first, we need a key syntactic lemma.

\begin{lemm}\label{syntactic}
If $\mathfrak w\in|\mathfrak Q(\Phi)|$ then 
\[\vdash\bigvee_{\mathfrak v\in\rho(\mathfrak w)}\mathrm{Sim}(\mathfrak v)\to\next\bigvee_{\mathfrak v\in\rho(\mathfrak w)}\mathrm{Sim}(\mathfrak v).\]
\end{lemm}

\proof
By Proposition \ref{propsub}.5 we have that, for all $\mathfrak v\in\rho(\mathfrak w)$,
\[\vdash\mathrm{Sim}(\mathfrak v)\to\next\bigvee_{\mathfrak v\dot\mapsto\mathfrak u}\mathrm{Sim}(\mathfrak u).\]

Now, if $\mathfrak v\dot\mapsto\mathfrak u$, it does not immediately follow that $\mathfrak u\in\rho(\mathfrak w)$, so we must examine the possible exceptions; let us see what these are.

Pick an efficient path $\vec{\mathfrak y}$ of length $N$ with $\mathfrak y_0=\mathfrak w$ and $\mathfrak y_{N-1}=\mathfrak v$.

Consider an extension $\vec{\mathfrak z}$ of $\vec{\mathfrak y}$ with $\mathfrak z_n=\mathfrak y_n$ for $n<N$ and $\mathfrak z_N=\mathfrak u$.

Now, by definition, the path $\vec{\mathfrak z}$ shows that $\mathfrak u$ is reachable from $\mathfrak w$, unless one of the following happens:
\begin{enumerate}
\item $\mathfrak u\not\in|\mathfrak Q(\Phi)|$;
\item the path $\vec{\mathfrak z}$ is inefficient.
\end{enumerate}

If 1 holds, then $\vdash\neg\mathrm{Sim}(\mathfrak u)$ and we can directly remove $\mathfrak u$ from the disjunction.

If 2 holds, there is $n<N$ such that $\mathfrak z_n\unlhd\mathfrak z_N$ (since $\vec{\mathfrak z}$ is inefficient but $\vec{\mathfrak y}$ is still efficient). But then, by Proposition \ref{propsub}.2, we have that
\[\vdash\mathrm{Sim}(\mathfrak u)\to\mathrm{Sim}(\mathfrak z_n),\]
and clearly $\mathfrak z_n$ is reachable from $\mathfrak w$. Therefore we can replace $\mathfrak u$ by $\mathfrak z_n$ in the disjunction.

We conclude that
\[\vdash\mathrm{Sim}(\mathfrak v)\to\next\bigvee_{\mathfrak u\in\rho(\mathfrak w)}\mathrm{Sim}(\mathfrak u),\]
as desired.

Since $\mathfrak v$ was arbitrary, this shows that
\[\vdash \bigvee_{\mathfrak v\in\rho(\mathfrak w)}\mathrm{Sim}(\mathfrak v)\to\next\bigvee_{\mathfrak v\in\rho(\mathfrak w)}\mathrm{Sim}(\mathfrak v).\]
\endproof

From this we obtain the following, which evidently implies $\omega$-sensibility:

\begin{prop}\label{tempinc}
If $\mathfrak w\in|\mathfrak Q(\Phi)|$ and $\ev\psi\in t(\mathfrak w)$, then there is $\mathfrak v\in\rho(\mathfrak w)$ such that $\psi\in t(\mathfrak v)$.
\end{prop}

\proof Towards a contradiction, assume that $\mathfrak w\in|\mathfrak Q(\Phi)|$ and $\ev\psi\in t(\mathfrak w)$ but, for all $\mathfrak v\in\rho(\mathfrak w)$, $\psi\not\in t(\mathfrak w)$.

By Lemma \ref{syntactic},
\[\vdash\bigvee_{\mathfrak v\in\rho(\mathfrak w)}\mathrm{Sim}(\mathfrak v)\to\next\bigvee_{\mathfrak v\in\rho(\mathfrak w)}\mathrm{Sim}(\mathfrak v).\]

But then we can use necessitation and $\mathsf {Ind}_\hf$ to show that
\[\vdash \bigvee_{\mathfrak v\in\rho(\mathfrak w)}\mathrm{Sim}(\mathfrak v)\to\hf\bigvee_{\mathfrak v\in\rho(\mathfrak w)}\mathrm{Sim}(\mathfrak v);\]
in particular,
\begin{equation}\label{other}
\vdash \mathrm{Sim}(\mathfrak w)\to\hf\bigvee_{\mathfrak v\in\rho(\mathfrak w)}\mathrm{Sim}(\mathfrak v).
\end{equation}




Now let $\mathfrak v\in\rho(\mathfrak w)$. By Proposition \ref{propsub}.1 and the assumption that $\psi\not\in t(\mathfrak v)$ we have that
\[\vdash \mathrm{Sim}(\mathfrak v)\to\neg\psi,\]
and since $\mathfrak v$ was arbitrary,
\[\vdash \bigvee_{\mathfrak v\in \rho(\mathfrak w)}\mathrm{Sim}(\mathfrak v)\to\neg\psi.\]
Using necessitation and distributibity we further have that
\[\vdash\hf \bigvee_{\mathfrak v\in \mathrm{Reach(\mathfrak w)}}\mathrm{Sim}(\mathfrak v)\to \hf\neg\psi.\]

This, along with (\ref{other}), shows that
\[\vdash\mathrm{Sim}(\mathfrak w)\to \hf\neg\psi;\]
however, once again by Proposition \ref{propsub}.1 and our assumption that $\ev\psi\in t(\mathfrak w)$ we have that $\vdash\mathrm{Sim}(\mathfrak w)\to\ev\psi$, which inconsistent with $\hf\neg\psi$, showing that $\vdash\neg\mathrm{Sim}(\mathfrak w)$.

But this contradicts the assumption that $\mathfrak w\in|\mathfrak Q(\Phi)|=\mathsf{Cons}(\Phi)$, and we conclude that there can be no such $\mathfrak w$.
\endproof

\begin{cor}\label{laststretch}
Given any finite set of formulas $\Phi$, $\mathfrak Q(\Phi)$ is a quasimodel.
\end{cor}

\proof
By Lemma \ref{niceprop}, $|\mathfrak Q(\Phi)|$ is open in $|\mathfrak I(\Phi)|$ and $\mapsto_{\mathfrak Q(\Phi)}$ is serial, while by Proposition \ref{tempinc}, $\mapsto_{\mathfrak Q(\Phi)}$ is $\omega$-sensible.

Thus $|\mathfrak Q(\Phi)|$ is regular, and it follows from Lemma \ref{regquasi} that $\mathfrak Q(\Phi)$ is a quasimodel.
\endproof

\section{Completeness}\label{secomp}

We are now ready to prove completeness of $\mathcal{DTL}^\ast$.

\begin{theorem}[Completeness]\label{theocomp}
If $\models\varphi$, then $\vdash\varphi$.
\end{theorem}

\proof
Suppose $\varphi$ is a consistent formula and let
\[W=\cbra \mathfrak w\in\mathfrak I_0(\varphi):\varphi\in t(\mathfrak w)\cket.\]

Then, by Proposition \ref{propsub}.4 we have that
\[\vdash\varphi\to\bigvee_{\mathfrak w\in W}\mathrm{Sim}(\mathfrak w).\]

Since $\varphi$ is consistent, it follows that some $\mathfrak w^\ast\in W$ is consistent and hence $\mathfrak w^\ast\in|\mathfrak Q(\varphi)|$. By Corollary \ref{laststretch}, $\mathfrak Q(\varphi)$ is a quasimodel, so that by Theorem \ref{second}, $\varphi$ is satisied in $\lim \mathfrak Q(\varphi)$.
\endproof

\section{Conclusions and future work}\label{conc}

The primary motivation for studying $\mathcal{DTL}$ with infinitary temporal modalities is to apply techniques from modal logic to the study of topological dynamics. The methods developed until now have allowed us to bring the model theory of modal logic and Kripke semantics into the study of such systems. We believe that the results presented in the current paper may turn the study of $\mathcal{DTL}$ in a new direction, where proof-theoretic methods take on a leading role.

However, there is much to be done. Indeed, dynamic topological systems appear in multiple branches of mathematics, but rarely is such a general class as that of {\em all} systems studied: applications typically consider systems with greater structure. Motivated by this fact, the semantic work in \cite{me2} has already been applied to metric spaces in \cite{me:metric} and minimal systems in \cite{me:minimal}. The advances presented here should lead to an analogous syntactic treatment of such classes of systems.

As a case in point, consider $\mathcal{DTL}$ over the class of minimal systems; a dynamic topological system $\mathfrak X$ is {\em minimal} if $|\mathfrak X|$ contains no non-empty, proper, closed, $f_\mathfrak X$-invariant subset. These systems have the property that the orbit of every point is dense in the whole space.

In \cite{me:minimal}, the language of $\mathcal{DTL}$ is enriched with a universal modality $\forall$. The formula
\[\exists\nc p\to\forall\ev p\]
can then be seen to be valid over this class; it expresses the fact that the orbit of every point is dense. Further, it is shown that

\begin{theorem}\label{mainmin}
Given a formula $\varphi$ in the (monadic) language $\mathsf L_{\ps\next\hf\forall}$, there exists a finite set of formulas $\mathrm{Ax}_\varphi$ of the form $\exists\nc\psi\to\forall\ev\psi$ such that $\varphi$ is valid over the class of all minimal dynamic topological systems if and only if $\bigwedge\mathrm{Ax}_\varphi\to\varphi$ is valid over the class of all dynamic topological systems.
\end{theorem}

This leads us to the following:

\begin{conj}
The set of $\mathsf L^\ast_{\ps\next\hf\forall}$-formulas valid over the class of minimal dynamic topological systems can be axiomatized by
\[\mathcal{DTL}^\ast+\mathsf{Ax}_\forall+\exists\nc p\to\forall\ev p,\]
where $\mathsf{Ax}_\forall$ denotes the standard axioms and rules for the universal modality and its interaction with other modalities.
\end{conj}

\proof[Proof idea.]
If a formula $\varphi$ is valid over the class of all minimal systems, it follows from Theorem \ref{mainmin} that $\bigwedge\mathrm{Ax}_\varphi\to\varphi$ is valid. Then, by Theorem \ref{theocomp},
\[\mathcal{DTL}^\ast\vdash\bigwedge\mathrm{Ax}_\varphi\to\varphi\]
and thus
\[\mathcal{DTL}^\ast+\exists\nc p\to\forall\ev p\vdash \varphi.\qedhere\]
\endproof

However, there are two obstacles to make this a proper proof:
\begin{enumerate}
\item Theorem \ref{theocomp} considers a language without the universal modality, and hence would have to be revised to accommodate for this addition.
\item Theorem \ref{mainmin} considers a monadic language, and would also have to be revised to include the polyadic $\ps$.
\end{enumerate}

It seems very unlikely that either of these points would prove terribly challenging: unfortunately, it also seems like the only way to deal with them (at least the second) would be to go back to the original proofs and check that the new additions do not pose a problem. Because of this, we shall not provide a full proof, but leave it as a conjecture and an indication of the new territory that must now be explored within the field of Dynamic Topological Logic.

\bibliographystyle{plain}
\bibliography{biblio}

\begin{thebibliography}{10}

\bibitem{arte}
S.N. Artemov, J.M. Davoren, and A.~Nerode.
\newblock Modal logics and topological semantics for hybrid systems.
\newblock {\em Technical Report MSI 97-05}, 1997.

\bibitem{do}
A.~Dawar and M.~Otto.
\newblock Modal characterisation theorems over special classes of frames.
\newblock {\em Annals of Pure and Applied Logic}, 161:1--42, 2009.
\newblock Extended journal version LICS 2005 paper.

\bibitem{me:minimal}
D.~Fern\'andez-Duque.
\newblock Dynamic topological logic interpreted over minimal systems.
\newblock {\em Journal of Philosophical Logic}.
\newblock Forthcoming.

\bibitem{me:simulability}
D.~Fern\'andez-Duque.
\newblock On the definability of simulability by finite $\mathsf{S4}$ models.
\newblock manuscript.

\bibitem{me2}
D.~Fern\'andez-Duque.
\newblock Non-deterministic semantics for dynamic topological logic.
\newblock {\em Annals of Pure and Applied Logic}, 157(2-3):110--121, 2009.
\newblock Kurt G\"odel Centenary Research Prize Fellowships.

\bibitem{me:metric}
D.~Fern\'andez-Duque.
\newblock Dynamic topological logic interpreted over metric spaces.
\newblock {\em Journal of Symbolic Logic}, 2011.

\bibitem{me:tangle}
D.~Fern\'andez-Duque.
\newblock Tangled modal logic for spatial reasoning, 2011.

\bibitem{dynamictangle}
D.~Fern\'andez-Duque.
\newblock Tangled modal logic for topological dynamics.
\newblock 2011.

\bibitem{pml}
D.~Gabelaia, A.~Kurucz, F.~Wolter, and M.~Zakharyaschev.
\newblock Non-primitive recursive decidability of products of modal logics with
  expanding domains.
\newblock {\em Annals of Pure and Applied Logic}, 142(1-3):245--268, 2006.

\bibitem{konev}
B.~Konev, R.~Kontchakov, F.~Wolter, and M.~Zakharyaschev.
\newblock Dynamic topological logics over spaces with continuous functions.
\newblock In G.~Governatori, I.~Hodkinson, and Y.~Venema, editors, {\em
  Advances in Modal Logic}, volume~6, pages 299--318, London, 2006. College
  Publications.

\bibitem{wolter}
B.~Konev, R.~Kontchakov, F.~Wolter, and M.~Zakharyaschev.
\newblock On dynamic topological and metric logics.
\newblock {\em Studia Logica}, 84:129--160, 2006.

\bibitem{s5}
P.~Kremer.
\newblock Dynamic topological $\mathsf{S5}$.
\newblock {\em Annals of Pure and Applied Logic}, 160:96--116, 2009.

\bibitem{kmints}
P.~Kremer and G.~Mints.
\newblock Dynamic topological logic.
\newblock {\em Annals of Pure and Applied Logic}, 131:133--158, 2005.

\bibitem{temporal}
O.~Lichtenstein and A.~Pnueli.
\newblock Propositional temporal logics: Decidability and completeness.
\newblock {\em Logic Jounal of the IGPL}, 8(1):55--85, 2000.

\bibitem{tarski}
A.~Tarski.
\newblock Der aussagenkalk\"ul und die topologie.
\newblock {\em Fundamenta Mathematica}, 31:103--134, 1938.

\end{thebibliography}

\end{document}